# $q$-Analogs for Steiner Systems and Covering Designs

Tuvi Etzion[*]    Alexander Vardy[†]

May 27, 2018

**Abstract**

The $q$-analogs of basic designs are discussed. It is proved that the existence of any unknown Steiner structures, the $q$-analogs of Steiner systems, implies the existence of unknown Steiner systems. Optimal $q$-analogs covering designs are presented. Some lower and upper bounds on the sizes of $q$-analogs covering designs are proved.

[*]Department of Computer Science, Technion, Haifa 32000, Israel, e-mail: etzion@cs.technion.ac.il.
[†]Department of Electrical Engineering, University of California San Diego, La Jolla, CA 92093, USA, e-mail: avardy@ucsd.edu.
[1]This research was supported in part by the United States — Israel Binational Science Foundation (BSF), Jerusalem, Israel, under Grant 2006097.

# 1 Introduction

Let $\mathbb{F}_q$ be a finite field with $q$ elements. For given integers $n \geq k \geq 0$, let $\mathcal{G}_q(n,k)$ denote the set of all $k$-dimensional subspaces of $\mathbb{F}_q^n$. $\mathcal{G}_q(n,k)$ is often referred to as Grassmannian. A *code* $\mathbb{C}$ over the grassmannian is a subset of $\mathcal{G}_q(n,k)$. In recent years there has been an increasing interest in codes over the Grassmannian as a result of their application to error-correction in network coding as was demonstrated by Koetter and Kschischang [10]. But, the interest in these codes has been also before this application, since these code are $q$-analogs of constant weight codes. Design theory is a well studied area in combinatorics related to coding theory. $q$-analogs of various combinatorial objects are well known. Many known combinatorial problems such as Sperner's Theorem [12] have $q$-analogs [20]. $q$-analogs of $t$-designs were studied in various papers and connections [2, 9, 13, 15, 16, 17, 18, 19]. These designs had also some interest in coding theory [1, 15]. In this paper we will consider $q$-analogs of covering designs for the first time to our knowledge. We will relate some of these coverings to other $q$-analogs of designs and codes known in the literature. We start our discussion with some preliminary definitions.

A *Steiner structure* $S_q[r,k,n]$ is a collection $\mathbb{S}$ of elements from $\mathcal{G}_q(n,k)$ such that each element from $\mathcal{G}_q(n,r)$ is contained in exactly one element of $\mathbb{S}$.

A *q-covering design* $C_q[n,k,r]$ is a collection $\mathbb{S}$ of elements from $\mathcal{G}_q(n,k)$ such that each element of $\mathcal{G}_q(n,r)$ is contained in at least one element of $\mathbb{S}$.

A *q-Turán design* $T_q[n,k,r]$ is a collection $\mathbb{S}$ of elements from $\mathcal{G}_q(n,r)$ such that each element of $\mathcal{G}_q(n,k)$ contains at least one element from $\mathbb{S}$.

The *q-covering number* $\mathcal{C}_q(n,k,r)$ is the minimum size of a $q$-covering design $C_q[n,k,r]$. The *q-Turán number* $\mathcal{T}_q(n,k,r)$ is the minimum size of a $q$-Turán design $T_q[n,k,r]$. Clearly, a Steiner structure $S_q[r,k,n]$ is the smallest $q$-covering design $C_q[n,k,r]$.

The goal of this paper is to examine these $q$-analogs combinatorial designs. The rest of the paper is organized as follows. In Section 2 we discuss the connections between $q$-covering designs and $q$-Turán designs. Basic lower bound on the $q$-covering numbers and its connection to Steiner structures is given. In Section 3 we prove a new necessary condition for the existence of Steiner structures. The condition implies that the task of constructing such structures with new parameters will be be very hard, to say the least. In Section 4 we present parameters for which we know the exact value of the $q$-covering numbers. In Section 5 we present a few lower bounds on the $q$-covering numbers which are better than the basic one. In Section 6 we present a recursive construction for $q$-covering designs which implies an upper bound on the $q$-covering numbers. In Section 7 we present constructions and bounds for specific values of $q$-covering numbers. Conclusion and problems for future research are given in Section 8

# 2 Relations between the Designs

In this section we will present simple connections between the various combinatorial objects which are discussed in this paper.

For a set $\mathbb{S} \subset \mathcal{G}_q(n,\ell)$ let $\mathbb{S}^\perp$, the *orthogonal complement* of $\mathbb{S}$, be the set

$$\mathbb{S}^\perp = \{A^\perp \ : \ A \in \mathcal{S}\},$$

where $A^\perp \in \mathcal{G}_q(n, n-\ell)$ is the *orthogonal complement* of the subspace $A$. For completeness we remark that all the results which follow hold for any $n$-dimensional subspace over $\mathbb{F}_q$ (and not just $\mathbb{F}_q^n$).

**Theorem 1.** $\mathbb{S}$ *is a $q$-covering design $C_q[n, k, r]$ if and only if $\mathbb{S}^\perp$ is a $q$-Turán design $T_q[n, n-r, n-k]$.*

*Proof.* Assume first that $\mathbb{S}$ is a $q$-covering design $C_q[n, k, r]$. Let $A \in \mathcal{G}_q(n, n-r)$; $A^\perp$ is an $r$-dimensional subspace and hence there exists at least one $k$-dimensional subspace $B \in \mathbb{S}$ such that $A^\perp \subset B$. This implies that $B^\perp \subset A$; since $B^\perp \in \mathbb{S}^\perp$ and $B^\perp \in \mathcal{G}_q(n, n-k)$, it follows that each subspace $A \in \mathcal{G}_q(n, n-r)$ contains at least one $(n-k)$-dimensional subspace of $\mathbb{S}^\perp$. Thus, $\mathbb{S}^\perp$ is a $q$-Turán design $T_q[n, n-r, n-k]$.

Similarly, if $\mathbb{S}$ is a $q$-Turán design $T_q[n, n-r, n-k]$ then $\mathbb{S}^\perp$ is a $q$-covering design $C_q[n, k, r]$. $\square$

**Corollary 1.** $\mathcal{C}_q(n, k, r) = \mathcal{T}_q(n, n-r, n-k)$.

Counting the number of number of subspaces in a $q$-covering designs involves the $q$-ary Gaussian coefficient $\begin{bmatrix} n \\ \ell \end{bmatrix}_q$ defined as follows (see [20]):

$$\begin{bmatrix} n \\ \ell \end{bmatrix}_q = \frac{(q^n - 1)(q^{n-1} - 1) \cdots (q^{n-\ell+1} - 1)}{(q^\ell - 1)(q^{\ell-1} - 1) \cdots (q - 1)}, \text{ and } \begin{bmatrix} n \\ 0 \end{bmatrix}_q = 1.$$

Now, we can prove a simple and basic lower bound on $\mathcal{C}_q(n, k, r)$.

**Lemma 1.** $\mathcal{C}_q(n, k, r) \geq \dfrac{\begin{bmatrix} n \\ r \end{bmatrix}_q}{\begin{bmatrix} k \\ r \end{bmatrix}_q}$ *with equality holds if and only if a Steiner structure $S_q[r, k, n]$ exists.*

*Proof.* Let $\mathbb{S}$ be a $q$-covering design $C_q[n, k, r]$ and let $P$ be a $k$-dimensional subspace of $\mathbb{S}$. $P$ contains $\begin{bmatrix} k \\ r \end{bmatrix}_q$ distinct $r$-dimensional subspaces. The $k$-dimensional subspaces of $\mathbb{S}$ must contain $\begin{bmatrix} n \\ r \end{bmatrix}_q$ distinct $r$-dimensional subspaces. Hence, $|\mathbb{S}| \geq \dfrac{\begin{bmatrix} n \\ r \end{bmatrix}_q}{\begin{bmatrix} k \\ r \end{bmatrix}_q}$. If $|\mathbb{S}| = \dfrac{\begin{bmatrix} n \\ r \end{bmatrix}_q}{\begin{bmatrix} k \\ r \end{bmatrix}_q}$ then each $r$-dimensional subspace is contained in exactly one element of $\mathbb{S}$, i.e., $\mathbb{S}$ is a Steiner structure $S_q[r, k, n]$. $\square$

## 3 On the Existence of Steiner Structure

A natural goal is to find optimal $q$-covering designs. From these designs the most interesting ones are the $q$-analogs of Steiner systems, i.e., the Steiner structure. For which parameters do Steiner structures exist? Clearly, $S_q[r, r, n]$ and $S_q[1, n, n]$ exist and these are trivial structures. The only known nontrivial Steiner structures are of the form $S_q[1, k, n]$, whenever $k$ divides $n$ [1, 3, 15]. These are also known as spreads and several constructions are known [3, 6, 8]. For all other parameters no Steiner structures are known and from computer searches we have performed and also ones reported by Thomas [19] it is tempting to conjecture that no such structures exists. There are many known results on Steiner structures [15, 19]. Two of them are given in the following two theorems.

**Theorem 2.** *If there exists a Steiner structure $S_q[t, k, n]$ then there exists a Steiner structure $S_q[t-1, k-1, n-1]$*

A *Steiner system* $S(r, k, n)$ is a collection $\mathcal{S}$ of $k$-subsets (called *blocks*) from an $n$-set such that each $r$-subset of the $n$-set is contained in exactly one block of $\mathcal{S}$.

**Theorem 3.**

- *If a Steiner structure $S_q[2, k, n]$ exists then there exists a Steiner system $S(2, \frac{q^k-1}{q-1}, \frac{q^n-1}{q-1})$.*

- *If a Steiner structure $S_q[2, k, n]$ exists then there exists a Steiner system $S(2, q^{k-1}, q^{n-1})$.*

- *If a Steiner structure $S_2[3, k, n]$ exists then there exists a Steiner system $S(3, 2^{k-1}, 2^{n-1})$.*

The connections between Steiner structures and Steiner systems given in Theorem 3 do not indicate the extreme difficulty to construct a Steiner structure $S_q[2, k, n]$. It is obvious that the most interesting and probably "easier to construct" case is $q = 2$. The following theorem will throw some light on the difficulty to construct such structures.

**Theorem 4.** *If a Steiner Structure $S_2[2, k, n]$ exists then a Steiner system $S(3, 2^k, 2^n)$ exists.*

*Proof.* Let $\mathbb{S}$ be a Steiner structure $S_2[2, k, n]$. We generate the following system of subsets over $\mathbb{F}_2^n$:

$$\mathcal{S} = \{\{\beta, \beta + \alpha_1, \beta + \alpha_2, \ldots, \beta + \alpha_{2^k-1}\} \ : \ \{0, \alpha_1, \ldots, \alpha_{2^k-1}\} \in \mathbb{S}, \ \beta \in \mathbb{F}_2^n\} \ .$$

First note that from each $k$-dimensional subspace of $\mathbb{S}$ we form $2^{n-k}$ disjoint $2^k$-subsets in $\mathcal{S}$. There are $\frac{\begin{bmatrix} n \\ 2 \end{bmatrix}_2}{\begin{bmatrix} k \\ 2 \end{bmatrix}_2}$ subspaces in $\mathbb{S}$ and therefore there are at most $\frac{2^{n-k}(2^n-1)(2^{n-1}-1)}{(2^k-1)(2^{k-1}-1)}$ blocks in $\mathcal{S}$. In a Steiner system $S(3, 2^k, 2^n)$ there are $\frac{\binom{2^n}{3}}{\binom{2^k}{3}} = \frac{2^{n-k}(2^n-1)(2^{n-1}-1)}{(2^k-1)(2^{k-1}-1)}$ blocks and hence to complete the proof we only have to show that each 3-subset $\{x, y, z\} \subset \mathbb{F}_2^n$ is contained in some block of $\mathcal{S}$. Since $\mathbb{S}$ is a Steiner structure $S_2[2, k, n]$, it follows that the two-dimensional subspace $\{0, x+y, x+z, y+z\}$ is contained in exactly one $k$-dimensional subspace $P$ of $\mathbb{S}$. Let $P = \{0, x+y, x+z, y+z, \gamma_1, \gamma_2, \ldots, \gamma_{2^k-4}\}$. By the definition of $\mathcal{S}$ we have that $P' = \{x, y, z, x+y+z, x+\gamma_1, x+\gamma_2, \ldots, x+\gamma_{2^k-4}\} \in \mathcal{S}$. Clearly $\{x, y, z\} \subset P'$ which completes the proof. $\square$

Theorem 4 is a strong evidence that constructing a Steiner structure $S_2[t,k,n]$ with new parameters is extremely difficult. If such a Steiner structure exists then by Theorem 2 there exists a Steiner structure $S_2[2, k-t+2, n-t+2]$. Hence by Theorem 4 there exists a Steiner system $S(3, 2^{k-t+2}, 2^{n-t+2})$, where $2^{k-t+2} \geq 8$. Unfortunately, no such Steiner systems are known, even so many efforts to find such systems were done for many years, which makes the task of finding new Steiner structures almost impossible. The Steiner structures which are the first target for investigation, due to their relative small parameters are $S_2[2,3,n]$. By the necessary conditions for their existence we must have $n \equiv 1$ or $3$ $(mod\ 6)$. Many other necessary conditions on the existence of such structures can be obtained by considering sets of derived designs related to any fixed $(n-1)$-dimensional subspace.

## 4 Optimal $q$-Covering Designs

In this section we we present a sequence of parameters for which we can compute the exact $q$-covering numbers. By Corollary 1, Lemma 1, and since a Steiner structure $S_q[1,k,n]$ exists if and only if $k$ divides $n$ we have the following theorem.

**Theorem 5.** $\mathcal{C}_q(n,k,1) = \mathcal{T}_q(n, n-1, n-k) = |S_q[1,k,n]| = \frac{q^n-1}{q^k-1}$, whenever $k$ divides $n$.

In the rest of this section we will find the exact value of $\mathcal{C}_q(n,k,1)$ and $\mathcal{C}_q(n, n-1, k)$ for all $q$, $n$, and $k$.

**Lemma 2.** $\mathcal{C}_q(n,k,1) \leq \mathcal{C}_q(n-1, k-1, 1)$.

*Proof.* We represent $\mathbb{F}_q^n$ by $\{(\alpha, \beta)\ :\ \alpha \in \mathbb{F}_q^{n-1},\ \beta \in \mathbb{F}_q\}$. Let $\mathbb{S}$ be a $q$-covering design $C_q[n-1, k-1, 1]$ and let $\mathbb{S}' = \{\{(\alpha, \beta)\ :\ \alpha \in P,\ \beta \in \mathbb{F}_q\}\ :\ P \in \mathbb{S}\}$. It is easy to verify that $\mathbb{S}'$ is a $q$-covering design $C_q[n,k,1]$ and the lemma follows. □

**Corollary 2.** If $\delta \geq 0$, then $\mathcal{C}_q(n+\delta, k+\delta, 1) \leq \mathcal{C}_q(n,k,1)$.

**Theorem 6.** If $2k \geq n$, then $\mathcal{C}_q(n,k,1) = q^{n-k} + 1 = \lceil \frac{q^n-1}{q^k-1} \rceil$.

*Proof.* By Theorem 5 we have that $\mathcal{C}(2(n-k), n-k, 1) = \frac{q^{2(n-k)}-1}{q^{n-k}-1} = q^{n-k} + 1$. Hence, if $2k \geq n$ then by Corollary 2 we have $\mathcal{C}_q(n,k,1) = \mathcal{C}(2(n-k) + 2k - n, n - k + 2k - n, 1) \leq \mathcal{C}(2(n-k), n-k, 1) = q^{n-k} + 1$.

By Lemma 1 we have that $\mathcal{C}_q(n,k,1) \geq \frac{\begin{bmatrix} n \\ 1 \end{bmatrix}_q}{\begin{bmatrix} k \\ 1 \end{bmatrix}_q} = \frac{q^n-1}{q^k-1}$ and hence $\mathcal{C}_q(n,k,1) \geq q^{n-k} + 1$.

Thus, $\mathcal{C}_q(n,k,1) = q^{n-k} + 1 = \lceil \frac{q^n-1}{q^k-1} \rceil$ if $2k \geq n$. □

**Corollary 3.** If $n \geq 2r$, then $\mathcal{T}_q(n, n-1, r) = q^r + 1$.

We continue to find the value of $\mathcal{C}_q(n,k,1)$ for $n > 2k$ by presenting a more general result on the size of an optimal $q$-covering design $C_q[n, \rho, 1]$ for $n$ which is not divisible by $\rho$. The first lemma is a generalization of constructions in [6, 7, 8] from $\mathbb{F}_2$ to $\mathbb{F}_q$. We omit the detailed proof which is almost identical to the one in [7].

**Lemma 3.** *Let $n = s\rho + m$, where $\rho < m < 2\rho$ and $s \geq 1$. There exist a set $\mathbb{S}$ with $\frac{q^{s\rho+m}-q^m}{q^\rho-1}$ $\rho$-dimensional subspaces of $\mathbb{F}_q^n$ and one $m$-dimensional subspace of $\mathbb{F}_q^n$ such that each one-dimensional subspace of $\mathbb{F}_q^n$ is a subspace of exactly one element from $\mathbb{S}$.*

**Corollary 4.** *If $1 \leq k \leq n$, then $\mathcal{C}_q(n, k, 1) = \lceil \frac{q^n-1}{q^k-1} \rceil$.*

*Proof.* If $k$ divides $n$ then the claim follows from Theorem 5. If $2k \geq n$ then the claim follows from Theorem 6. Hence we only have to consider the case $n = sk + m$, where $k < m < 2k$ and $s \geq 1$. By Lemma 3 we have that

$$\mathcal{C}(sk+m, k, 1) \leq \frac{q^{sk+m} - q^m}{q^k - 1} + \mathcal{C}(m, k, 1) \tag{1}$$

and by applying Theorem 6 in (1) we obtain

$$\mathcal{C}(sk+m, k, 1) \leq q^{(s-1)k+m} + q^{(s-2)k+m} + \cdots + q^m + q^{m-k} + 1 = \left\lceil \frac{q^{sk+m} - 1}{q^k - 1} \right\rceil.$$

□

**Corollary 5.** *If $1 \leq r \leq n - 1$, then $\mathcal{T}_q(n, n-1, r) = \lceil \frac{q^n-1}{q^{n-r}-1} \rceil$.*

We now know the value of $\mathcal{C}_q(n, k, 1)$ for any given $k$, $1 \leq k \leq n$. We will consider now the "complement" value $\mathcal{C}_q(n, n-1, r)$ for any given $1 \leq r \leq n-1$. By Corollary 1 we have $\mathcal{C}_q(n, n-1, r) = \mathcal{T}_q(n, n-r, 1)$ and hence we to consider the value $\mathcal{T}_q(n, k, 1)$ for any given $1 \leq k \leq n$.

We start by presenting a basic upper bound on on $\mathcal{T}_q(n, k, r)$.

**Lemma 4.** $\mathcal{T}_q(n, k, r) \leq \begin{bmatrix} n - k + r \\ r \end{bmatrix}_q.$

*Proof.* Let $Q$ be any $(n-k+r)$-dimensional subspace of $\mathbb{F}_q^n$ and let $\mathbb{S}$ consists of all the $r$-dimensional subspaces of $Q$, i.e., $|\mathbb{S}| = \begin{bmatrix} n-k+r \\ r \end{bmatrix}_q$. Let $P$ be a $k$-dimensional subspace of $\mathbb{F}_q^n$. Since $\dim Q + \dim P = n + r$ it follows that $\dim(P \cap Q) \geq r$. Hence, $P$ contains an element from $\mathbb{S}$. Thus, $\mathbb{S}$ is a $q$-Turán design $T_q[n, k, r]$ and $\mathcal{T}_q(n, k, r) \leq \begin{bmatrix} n-k+r \\ r \end{bmatrix}_q.$ □

**Corollary 6.** $\mathcal{C}_q(n, k, r) = \mathcal{T}_q(n, n-r, n-k) \leq \begin{bmatrix} n - (n-r) + (n-k) \\ n-k \end{bmatrix}_q = \begin{bmatrix} n-k+r \\ r \end{bmatrix}_q.$

In the sequel let $\langle A \rangle$ denote the subspace of $\mathbb{F}_q^n$ spanned by the elements of a set $A \subseteq \mathbb{F}_q^n$.

**Theorem 7.** *If $1 \leq k \leq n$, then $\mathcal{T}_q(n, k, 1) = \frac{q^{n-k+1}-1}{q-1}$.*

*Proof.* By Lemma 4 we have $\mathcal{T}_q(n, k, 1) \leq \begin{bmatrix} n-k+1 \\ 1 \end{bmatrix}_q = \frac{q^{n-k+1}-1}{q-1}$.

Now, assume that $\mathbb{S}$ is a set with $\frac{q^{n-k+1}-1}{q-1} - 1$ one-dimensional subspaces of $\mathbb{F}_q^n$. We will show that the exists a $k$-dimensional subspace of $\mathbb{F}_q^n$ which does not contain any subspace of $\mathbb{S}$. Let $Q$ be the largest subspace of $\mathbb{F}_q^n$ which does not contain any subspace of $\mathbb{S}$, $\ell = \dim Q$. If $\ell \geq k$ then any $k$-dimensional subspace of $Q$ does not contain any subspace of $\mathbb{S}$ and our claim is proved. Hence, we assume for contradiction that $\ell < k$. $Q$ consists of $q^\ell - 1$ nonzero elements. Let $B \in \mathbb{S}$; $\langle B \cup Q \rangle$ has $q^{\ell+1} - q^\ell$ nonzero elements which are not contained in $Q$. Hence, $|\cup_{B \in \mathbb{S}} \langle B \cup Q \rangle| \leq (\frac{q^{n-k+1}-1}{q-1} - 1)(q^{\ell+1} - q^\ell) + q^\ell - 1 = (q^{n-k+1} - q)q^\ell + q^\ell - 1 = q^{n-k+\ell+1} - q^{\ell+1} + q^\ell - 1$. Since $\ell < k$ it follows that $q^{n-k+\ell+1} - q^{\ell+1} + q^\ell - 1 < q^n - 1$ and hence there exist a nonzero element $\alpha \in \mathbb{F}_q^n$ such that $\alpha \notin \cup_{B \in \mathbb{S}} \langle B \cup Q \rangle$. Therefore $\langle \{\alpha\} \cup Q \rangle$ does not contain a subspace from $\mathbb{S}$ and its dimension is $\ell + 1$ in contradiction to the maximality of $Q$. Hence, there exists a $k$-dimensional subspace of $\mathbb{F}_q^n$ which does not contain any subspace of $\mathbb{S}$ and therefore, $\mathcal{T}_q(n, k, 1) \geq \frac{q^{n-k+1}-1}{q-1}$.

Thus, $\mathcal{T}_q(n, k, 1) = \frac{q^{n-k+1}-1}{q-1}$. □

**Corollary 7.** $\mathcal{C}_q(n, n-1, k) = \frac{q^{k+1}-1}{q-1}$.

## 5 Lower Bounds on the $q$-Covering Numbers

In this section we present lower bounds on the $q$-covering numbers. The most basic bound is given in Lemma 1. We will derive now a bound for which the basic bound is a special case. We start by proving a $q$-analog for the known Schönheim bound [14].

**Theorem 8.** $\mathcal{C}_q(n, k, r) \geq \lceil \frac{q^n-1}{q^k-1} \mathcal{C}_q(n-1, k-1, r-1) \rceil$.

*Proof.* Let $\mathbb{S}$ be a $q$-covering design which attains the value of $\mathcal{C}_q(n, k, r)$. Each element of $\mathbb{S}$ contains $\frac{q^k-1}{q-1}$ one-dimensional subspaces of $\mathbb{F}_q^n$. There are $\frac{q^n-1}{q-1}$ one-dimensional subspaces in $\mathbb{F}_q^n$. Hence, there exists an one-dimensional subspace $P$ of $\mathbb{F}_q^n$ which is contained in at most $\frac{q^k-1}{q^n-1} \mathcal{C}_q(n, k, r)$ elements of $\mathbb{S}$. Assume $P$ is contained in $\ell$ elements of $\mathbb{S}$. Let $\mathbb{F}_q^n = P \oplus Q$, where $Q$ is an $(n-1)$-dimensional subspace of $\mathbb{F}_q^n$. We define the following set $\mathbb{S}'$.

$$\mathbb{S}' = \{Y \cap Q \,:\, Y \in \mathbb{S},\, P \in Y\}.$$

Clearly $\mathbb{S}'$ contains $\ell$ $(k-1)$-dimensional subspaces of $Q$, one subspace for each element of $\mathbb{S}$ which contains $P$.

Let $A$ be an $(r-1)$-dimensional subspace of $Q$. $P \oplus A$ is an $r$-dimensional subspace of $\mathbb{F}_q^n$ and hence there exists a element $X \in \mathbb{S}$ such that $P \oplus A \subset X$. Clearly, $X \cap Q \in \mathbb{S}'$ and $A = (P \oplus A) \cap Q \subset X \cap Q$, i.e., $\mathbb{S}'$ is a $q$-covering design $\mathcal{C}_q[n-1, k-1, r-1]$. Therefore we have

$$\frac{q^k-1}{q^n-1} \mathcal{C}_q(n, k, r) \geq \mathcal{C}_q(n-1, k-1, r-1),$$

which implies

$$\mathcal{C}_q(n, k, r) \geq \frac{q^n-1}{q^k-1} \mathcal{C}_q(n-1, k-1, r-1) \,.$$

□

**Corollary 8.** $\mathcal{C}_q(n,k,r) \geq \lceil \frac{q^n-1}{q^k-1} \lceil \frac{q^{n-1}-1}{q^{k-1}-1} \cdots \lceil \frac{q^{n+1-r}-1}{q^{k+1-r}-1} \rceil \cdots \rceil \rceil \geq \frac{\begin{bmatrix} n \\ r \end{bmatrix}_q}{\begin{bmatrix} k \\ r \end{bmatrix}_q}$.

*Proof.* The result is an immediate consequence from iterative applications of Theorem 8 and Corollary 4. □

It is interesting to note that Lemma 1 is a special case of Corollary 8.

We will now modify a bound given in [5] which is a special case of a theorem in [4]. The following theorem is a q-analog of a related theorem in [5]. The proof is q-identical to the proof in [5] and hence it will be omitted. The theorem is stated as in [5] in terms of Turán numbers.

**Theorem 9.** $\mathcal{T}_q(n, r+1, r) \geq \frac{(q^{n-r}-1)(q-1)}{(q^r-1)^2} \begin{bmatrix} n \\ r-1 \end{bmatrix}_q$.

**Corollary 9.** $\mathcal{C}_q(n, k, k-1) \geq \frac{(q^k-1)(q-1)}{(q^{n-k}-1)^2} \begin{bmatrix} n \\ k+1 \end{bmatrix}_q$.

## 6 An Upper Bound on the $q$-Covering Numbers

In this section we present a simple construction which imply a relatively good upper bound on the $q$-covering numbers. Although the bound of Lemma 4 (Corollary 6) is attained (see Theorem 7 and Corollary 7) it is usually a weak bound. In the sequel this bound is considerably improved.

**Theorem 10.** $\mathcal{C}_q(n,k,r) \leq q^{n-k}\mathcal{C}_q(n-1,k-1,r-1) + \mathcal{C}_q(n-1,k,r)$.

*Proof.* We represent $\mathbb{F}_q^n$ by $\{(\alpha, \beta) : \alpha \in \mathbb{F}_q^{n-1}, \beta \in \mathbb{F}_q\}$. Let $\mathbb{S}_1$ be a $q$-covering design $C_q[n-1, k-1, r-1]$ and $\mathbb{S}_2$ be a $q$-covering design $C_q[n-1, k-1, r-1]$. We form a set $\mathbb{S}$ from the following two types of $k$-dimensional subspaces of $\mathbb{F}_q^n$.

[V.1] For each subspace $P = \{0, \alpha_1, \cdots, \alpha_{q^{k-1}-1}\} \in \mathbb{S}_1$ let $P_1 = P, P_2, \ldots, P_{q^{n-k}}$ be the disjoint cosets of $P$ in $\mathbb{F}_q^{n-1}$. Let $\beta_0 = 0, \beta_1, \ldots, \beta_{q^{n-k}}$ be any $q^{n-k}$ coset representatives, i.e., $\beta_i \in P_i$, $1 \leq i \leq q^{n-k}$. For each $1 \leq i \leq q^{n-k}$ we form the subspace $\langle\{(\alpha_1, 0), \cdots, (\alpha_{q^{k-1}-1}, 0), (\beta_i, 1)\}\rangle$ in $\mathbb{S}$.

[V.2] For each subspace $\{0, \alpha_1, \cdots, \alpha_{q^k-1}\} \in \mathbb{S}_2$ the subspace $\{(0,0), (\alpha_1, 0), \cdots, (\alpha_{q^k-1}, 0)\}$ is formed in $\mathbb{S}$.

It is clear that that $\mathbb{S}$ consists of $q^{n-k}|\mathbb{S}_1| + |\mathbb{S}_2|$ $k$-dimensional subspaces. Therefore, to complete the proof we only have to show that for any given $r$-dimensional subspace $T$ of $\mathbb{F}_q^n$ there exists a subspace $Q \in \mathbb{S}$ such that $T \subset Q$. We distinguish between two cases:
**Case 1:** If $T \subset \mathbb{F}_q^{n-1} \times \{0\}$ then clearly $T$ is contained in a subspace formed by [V.2] since $\mathbb{S}_2$ is a $q$-covering design $C_q[n-1, k, r]$.
**Case 2:** If $T = \langle\{(\alpha_1, 0), \cdots, (\alpha_{q^{r-1}-1}, 0), (\beta, \gamma)\}\rangle$, where $\beta \in \mathbb{F}_q^{n-1}$ then w.l.o.g. we can assume that $\gamma = 1$. Consider the $r$-dimensional subspace $T' = \{0, \alpha_1, \cdots, \alpha_{q^{r-1}-1}\}$. Since $\mathbb{S}_1$ is a $q$-covering design $C_q[n-1, k-1, r-1]$ it follows that there exists a $(k-1)$-dimensional subspace $Q' \in \mathbb{S}_1$ such that $T' \subset Q'$. Let $Q = \langle \{(\alpha, 0) : \alpha \in Q'\} \cup \{(\beta, 1)\} \rangle$. Clearly $Q$ is a $k$-dimensional subspace such that $Q \in \mathbb{S}$ by [V.2] and $T \subset Q$.

Cases 1 and 2 implies that $\mathbb{S}$ is a $q$-covering design $C_q[n, k, r]$ and the theorem follows. $\square$

The construction of Theorem 10 can be used to obtain upper bounds on the covering numbers for many sets of parameters. For example, consider the recursion $g(n) = 4g(n-1) + 2^{n-2} - 1$, where $g(4) = 5$. The solution for this recursion is $g(n) = 9 \cdot 2^{2n-8} - 2^{n-2} - \frac{2^{2n-8}-1}{3}$. By Theorem 5 we have $\mathcal{C}_2(4, 2, 1) = 5$ and hence by Theorem 10 we have $\mathcal{C}(n, n-2, n-3) \leq g(n)$. By using better initial conditions for $g(n)$ as the ones obtained in the next section, this bound can be improved.

## 7 Covering Numbers for Small Parameters

In this section we show how some lower and upper bounds of Sections 5 and 6 can be improved. We will consider only the case where $q = 2$. We start with a specific lower bound which implies many other new lower bounds by using Theorem 8.

**Lemma 5.** $\mathcal{C}_2(5, 3, 2) \geq 27$.

*Proof.* Let $\mathbb{S}$ be a $q$-covering design $C_2[5, 3, 2]$ which attains $\mathcal{C}_2(5, 3, 2)$. For each $x \in \mathbb{F}_2^5 \setminus \{0\}$ let $\eta(x) = |\{P : x \in P, P \in \mathbb{S}\}|$. For each $x, y \in \mathbb{F}_2^5 \setminus \{0\}$ there exists a 2-dimensional subspace of $\mathbb{F}_2^5$ which contains $x$ and $y$. Hence, since for each $P \in \mathbb{S}$, $x \in P$, we have $|P \setminus \{0, x\}| = 6$, it follows that $\eta(x) \geq \frac{30}{6} = 5$ for each $x \in \mathbb{F}_2^5 \setminus \{0\}$. We distinguish between two cases:

**Case 1:** For each $x \in \mathbb{F}_2^5 \setminus \{0\}$, $\eta(x) \geq 6$. Hence, $\sum_{x \in \mathbb{F}_2^5 \setminus \{0\}} \eta(x) \geq 31 \cdot 6 = 186$ and hence $\mathcal{C}_2(5, 3, 2) \geq \lceil \frac{186}{7} \rceil = 27$.

**Case 2:** Let $\eta(z) = 5$ for some $z \in \mathbb{F}_2^5 \setminus \{0\}$ and let $Q_i$, $1 \leq i \leq 5$, be a subset such that $|Q_i| = 6$ and $Q_i \cup \{0, z\} \in \mathbb{S}$. It follows that for each $i \neq j$ we have $|Q_i \cap Q_j| = 0$ and $\bigcup_{i=1}^{5} Q_i = \mathbb{F}_2^5 \setminus \{0, z\}$.

**Fact A:** If $P \in \mathbb{S}$, $P \notin \{Q_1, Q_2, Q_3, Q_4, Q_5\}$, then by the pigeon hole principle we have that there exists a $j$ such that $|Q_j \cap P| = 3$. Since $\Sigma_{i=0}^{5}|Q_i \cap P| = 7$; for each $i$, $|Q_i \cap P| > 0$; and for each $i \neq \ell$, $|Q_i \cap Q_\ell| = 0$, it follows for $\ell \neq j$ that $|Q_\ell \cap P| = 1$.

Let $\sigma(Q_i) = |\{P \in \mathbb{S} : |Q_i \cap P| = 3\}|$, $1 \leq i \leq 5$, and w.l.o.g. we assume that $\sigma(Q_1) \leq \sigma(Q_i)$, $2 \leq i \leq 5$. The following three claims will lead to the proof of the lemma.

**Claim 1:** If $P \in \mathbb{S}$, $|Q_1 \cap P| = 3$, and $x \in Q_1 \cap P$ then $|\{z\} \cup Q_1 \cup P \setminus \{0\}| = 11$ and hence $\eta(x) \geq 2 + \lceil \frac{20}{6} \rceil = 6$.

**Claim 2:** If $P_1, P_2 \in \mathbb{S}$, $|Q_1 \cap P_1| = |Q_1 \cap P_2| = 3$, and $x \in Q_1 \cap P_1 \cap P_2$ then $\eta(x) \geq 7$.

Assume the contrary, that $\eta(x) < 7$; therefore, by Claim 1 we have that $\eta(x) = 6$. W.l.o.g. we assume that $P_3, P_4, P_5 \in \mathbb{S}$ and $x \in P_i$, $3 \leq i \leq 5$. By Fact A, for each $P_i$, $1 \leq i \leq 5$, there exists exactly one $Q_j$, $1 \leq j \leq 5$, such that $|P_i \cap Q_j| = 3$. Therefore, for one of $Q_2, Q_3, Q_4, Q_5$, w.l.o.g. $Q_2$, we have that $|Q_2 \cap P_i| = 1$ for $1 \leq i \leq 5$. Hence, there exists a $y \in Q_2$, such that $y \notin Q_1 \cup \bigcup_{i=1}^{5} P_i$. It implies that there exists $P_6 \in \mathbb{S}$ such that

$\{x, y\} \subset P_6$. Thus, $\eta(x) \geq 7$.

**Claim 3:** If $2 \leq \sigma(Q_1) \leq 4$ then $\mathcal{C}_2(5, 3, 2) \geq 35 - 2\sigma(Q_1)$.

1. Let $x \in Q_1$ and $\lambda(x) = |\{P \in \mathbb{S} : |Q_1 \cap P| = 3, \ x \in P\}|$. Since $\sigma(Q_1) \leq 4$ one can easily verify that in the set $A = \{x : \lambda(x) \leq 2\}$ ($A$ includes the elements of $Q_1$ which belong to at most two subspaces of $\mathbb{S}$ intersecting $Q_1$ in three elements) there are at least 3 elements. For $y \in A$ we consider the value $\eta(y) - \lambda(y)$. By Claims 1 and 2 and since $\eta(y) \geq 5$, it follows that this value is at least 5, i.e., $|\{P \in \mathbb{S} : Q_1 \cap P = \{y\}\}| = \eta(y) - \lambda(y) - 1 \geq 4$.

2. $2 \leq \sigma(Q_1)$ implies that $6 - \sigma(Q_1) \leq 4$ and by Claim 2, for $y \in Q_1 \setminus A$ we have $\eta(y) \geq 7$ and $|\{P \in \mathbb{S} : Q_1 \cap P = \{y\}\}| = \eta(y) - \lambda(y) - 1 \geq 6 - \sigma(Q_1)$.

Therefore, we have $\mathcal{C}_2(5, 3, 2) \geq \eta(z) + \sigma(Q_1) + 3 \cdot 4 + 3 \cdot (6 - \sigma(Q_1))$ and the claim follows.

After establishing the three claims we continue by considering the value of $\sigma(Q_1)$.

- If $\sigma(Q_1) = 0$ then $\mathcal{C}_2(5, 3, 2) = |\mathbb{S}| \geq \eta(z) + \sum_{x \in Q_1}(\eta(x) - 1) \geq 5 + 6 \cdot 4 = 29$.

- If $\sigma(Q_1) = 1$ then $\mathcal{C}_2(5, 3, 2) = |\mathbb{S}| \geq \eta(z) + \sum_{x \in Q_1}(\eta(x) - 1) - 2 \geq 5 + 6 \cdot 4 - 2 = 27$.

- If $2 \leq \sigma(Q_1) \leq 4$ then by Claim 3 we have $\mathcal{C}_2(5, 3, 2) \geq 27$.

- If $\sigma(Q_1) \geq 5$ then $\sigma(Q_i) \geq 5$, $2 \leq i \leq 5$. By Fact A each $P \in \mathbb{S}$ intersect exactly one $Q_i$ in three elements. Hence, we have $\mathcal{C}_2(5, 3, 2) = |\mathbb{S}| \geq \eta(z) + \sum_{i=1}^{5} \sigma(Q_i) \geq 5 + 5 \cdot \sigma(Q_1) \geq 30$.

Thus, $\mathcal{C}_2(5, 3, 2) \geq 27$. □

As a consequence of Lemma 5, Theorem 10, Corollary 7, and Theorem 5 we have.

**Corollary 10.** $\mathcal{C}_2(5, 3, 2) = 27$.

Now, we improve one important upper bound. By using Theorem 10 and this bound many more new upper bounds can be obtained.

**Lemma 6.** $\mathcal{C}_2(7, 3, 2) \leq 399$.

*Proof.* Let $\alpha$ be a root of the primitive polynomial $x^6 + x + 1$, and hence a primitive element in GF(64). $\mathbb{F}_2^7$ will be represented as $\mathbb{F}_2^7 = \{(\beta, 0) : \beta \in GF(64)\} \cup \{(\beta, 1) : \beta \in GF(64)\}$.

Consider the five subsets $A_1 = \{0, 1, 4, 16\}$, $A_2 = \{0, 2, 8, 32\}$, $A_3 = \{0, 5, 27, 40\}$, $A_4 = \{0, 7, 44, 53\}$, $A_5 = \{0, 11, 29, 49\}$ and let $D_i = \{a - b : a, b \in A_i, a \neq b\}$ for $1 \leq i \leq 5$. Clearly, $|D_i| = 12$, $1 \leq i \leq 5$, and $\bigcup_{i=1}^{5} D_i = \mathbb{Z}_{63} \setminus \{0, 21, 42\}$. Moreover, if $A_i = \{j_1, j_2, j_3, j_4\}$, $1 \leq i \leq 5$, then $\alpha^{j_1} + \alpha^{j_2} + \alpha^{j_3} + \alpha^{j_4} = 0$. Now, we form a set $\mathcal{A}$ with 336 three-dimensional subspaces as follows.

For each $i$, $1 \leq i \leq 5$, $A_i = \{j_1, j_2, j_3, j_4\}$, and $\ell$, $0 \leq \ell \leq 62$, we form in $\mathcal{A}$ 63 three-dimensional subspaces of the form

$$\{(0,0), (\alpha^{j_1+\ell} + \alpha^{j_2+\ell}, 0), (\alpha^{j_1+\ell} + \alpha^{j_3+\ell}, 0), (\alpha^{j_1+\ell} + \alpha^{j_4+\ell}, 0), (\alpha^{j_1+\ell}, 1), (\alpha^{j_2+\ell}, 1), (\alpha^{j_3+\ell}, 1), (\alpha^{j_4+\ell}, 1)\}.$$

For each $\ell$, $0 \leq \ell \leq 20$ we form in $\mathcal{A}$ the 21 three-dimensional subspaces of the form

$$\{(0,0), (\alpha^\ell, 0), (\alpha^{21+\ell}, 0), (\alpha^{42+\ell}, 0), (0,1), (\alpha^\ell, 1), (\alpha^{21+\ell}, 1), (\alpha^{42+\ell}, 1)\} \ .$$

Clearly, each pair $((\alpha^{j_1}, 1), (\alpha^{j_2}, 1))$ is contained in exactly one subspace of $\mathcal{A}$, and hence each two-dimensional subspace $\{(0,0), (\alpha^{j_1}, 1), (\alpha^{j_2}, 1), (\alpha^{j_1} + \alpha^{j_2}, 0)\}$ is a subspace of exactly one three-dimensional subspace of $\mathcal{A}$. Similarly, each two-dimensional subspace $\{(0,0), (0,1), (\alpha^i, 1), (\alpha^i, 0)\}$ is a subspace of exactly one three-dimensional subspace of $\mathcal{A}$.

Finally, each one of the 273 two-dimensional subspaces of the form $\{(0,0), (\alpha^\ell, 0), (\alpha^{\ell+7}, 0), (\alpha^{\ell+26}, 0)\}$, $\{(0,0), (\alpha^\ell, 0), (\alpha^{\ell+11}, 0), (\alpha^{\ell+25}, 0)\}$, $\{(0,0), (\alpha^\ell, 0), (\alpha^{\ell+13}, 0), (\alpha^{\ell+35}, 0)\}$, $\{(0,0), (\alpha^\ell, 0), (\alpha^{\ell+9}, 0), (\alpha^{\ell+45}, 0)\}$, $0 \leq \ell \leq 62$, and $\{(0,0), (\alpha^\ell, 0), (\alpha^{\ell+21}, 0), (\alpha^{\ell+42}, 0)\}$, $0 \leq \ell \leq 20$, is contained in one three-dimensional subspace of $\mathcal{A}$.

To the set $\mathcal{A}$ we add a set $\mathcal{B}$ with 63 three-dimensional subspaces of the form $\{(0,0), (\alpha^\ell, 0), (\alpha^{1+\ell}, 0), (\alpha^{4+\ell}, 0), (\alpha^{6+\ell}, 0), (\alpha^{16+\ell}, 0), (\alpha^{24+\ell}, 0), (\alpha^{33+\ell}, 0)\}$, $0 \leq \ell \leq 62$. These 63 three-dimensional subspaces contain the 378 two-dimensional subspaces of the form $\{(0,0), (\alpha^\ell, 0), (\alpha^{2^i+\ell}, 0), (\alpha^{3 \cdot 2^{i+1}+\ell}, 0)\}$, $0 \leq i \leq 5$, $0 \leq \ell \leq 62$.

Therefore, each one of the $\begin{bmatrix} 7 \\ 2 \end{bmatrix} = \begin{bmatrix} 6 \\ 2 \end{bmatrix} + \binom{64}{2} = 3 \cdot 7 \cdot 127$ two-dimensional subspaces of $\mathbb{F}_2^7$ is contained in one of the 399 three-dimensional subspace. Thus, $\mathcal{C}_2(7,3,2) \leq 399$. □

## 8 Conclusion and Problems for Future Research

The $q$-analog of covering designs and the $q$-covering number $\mathcal{C}_q(n,k,r)$ were defined. The covering numbers were completely determined for $r = 1$ and for $k = n - 1$. Basic lower and upper bounds these numbers were derived. The existence of Steiner structures were considered and a strong necessary condition for their existence was proved. We conclude with a list of three problems for further research.

- Does there exists a Steiner structure $S_2[2, 3, 7]$? Thomas [18, 19] has examined this structure and didn't come to any conclusion concerning its existence. The covering design $C_2[7, 3, 2]$ of Size 399 constructed in Section 7 is close to such structure. This structure if exists will have 381 three-dimensional subspaces. A related packing of size 304 was found [11]. Improving on the size of this packing is also a step in this direction.

- The next set of values for which one might expect to to find the exact $q$-covering number is $\mathcal{C}_2(n, n-2, 2)$. By Theorem 8 and Corollary 4, a lower bound on this number, for $n \geq 6$, is $\mathcal{C}_2(n, n-2, 2) \geq \lceil \frac{2^n-1}{2^{n-2}-1} \mathcal{C}_2(n-1, n-3, 1) \rceil = \lceil \frac{2^n-1}{2^{n-2}-1} \lceil \frac{2^{n-1}-1}{2^{n-3}-1} \rceil \rceil = 21$. By Theorem 10 an upper bound on this number, for $n \geq 6$, is $\mathcal{C}_2(n, n-2, 2) \leq 4\mathcal{C}_2(n-1, n-3, 1) + \mathcal{C}_2(n-1, n-2, 2) = \lceil \frac{2^{n-1}-1}{2^{n-3}-1} \rceil + 7 = 27$. For $n = 5$ we have by Lemma 5 that $\mathcal{C}_2(5, 3, 2) = 27$. Can the lower bound be improved to get the same result for general $n$?

- The construction method described in Lemma 6 employs cyclic shifts on sets of elements in GF(64), where the elements on which the covering is preformed are in $\mathbb{F}_2^7$. Codes for packing which are cyclic were discussed in [7, 11]. Is there a general construction method for cyclic $q$-covering designs?